\documentclass{gtart}

\input gtoutput
\volumenumber{4}\papernumber{5}\volumeyear{2000}
\pagenumbers{171}{178}\published{19 July 2000}
\proposed{Walter Neumann}\seconded{Jean-Pierre Otal, Robion Kirby}
\received{18 June 2000}
\accepted{19 July 2000}

\newtheorem{theorem}{Theorem}[section]    
         
\newtheorem{corollary}[theorem]{Corollary}

 \def\SO{\mbox{\rm{SO}}}

 \def\A{\cal A} 
\def\sign{\mbox{\rm{sign}}}  
\def\correction{\mbox{\rm{correction}}} 
\def\Vol{\mbox{\rm{Vol}}}

\title{On the geometric boundaries of hyperbolic\\$4$--manifolds} 
\asciititle{On the geometric boundaries of hyperbolic 4-manifolds} 
\shorttitle{Geometric boundaries of hyperbolic 4--manifolds} 
\author{D\thinspace D Long\\A\thinspace W Reid}
\asciiauthors{DD Long\\AW Reid}
\address{Department of Mathematics, University of California\\
Santa Barbara, CA 93106, USA}
\secondaddress{Department of Mathematics, University of Texas\\ 
Austin, TX 78712, USA}
\email{long@math.ucsb.edu\\areid@math.utexas.edu}
\asciiaddress{Department of Mathematics, University of California\\
Santa Barbara, CA 93106, USA\\
Department of Mathematics, University of Texas\\ 
Austin, TX 78712, USA}

\begin{abstract}
We provide, for hyperbolic and flat 3--manifolds, obstructions to
bounding hyperbolic 4--manifolds, thus resolving in the negative a
question of Farrell and Zdravkovska.
\end{abstract}

\asciiabstract{We provide, for hyperbolic and flat 3-manifolds,
obstructions to bounding hyperbolic 4-manifolds, thus resolving in the
negative a question of Farrell and Zdravkovska.}




\begin{document} \maketitlepage 

\section{Introduction}

It is a classical result of Rohlin that the bordism group of closed
orientable $3$--manifolds is zero, so that every such $M^3$ can be
identified with $\partial W^4$ for some appropriate compact $W^4$.

This paper deals with two geometric incarnations of this
situation. The first, motivated in part by considerations in physics, (see
\cite{Gi}, \cite{RT} and \cite{W}), asks whether a closed orientable
hyperbolic $3$--manifold can be the totally geodesic boundary of a
compact complete hyperbolic $4$--manifold. The second, motivated by a
theorem of Hamrick and Royster \cite{HR} (which shows that that in
every dimension a {\em flat} $k$--manifold is nullbordant), concerns
the question of whether every flat $3$--manifold is, up to homeomorphism, a cusp
cross-section of a complete finite volume 1--cusped hyperbolic
$4$--manifold. 

We introduce the following notation.
If a hyperbolic 3--manifold $M$ is the totally geodesic boundary of a
hyperbolic $4$--manifold $W$, we say $M$ {\em bounds geometrically}. 
Also recall that the $\eta$--invariant of $M$, denoted $\eta (M)$
is defined as  $\eta(0)$
where $\eta(s)$ is formed from a signed collection of
eigenvalues of a certain first order self adjoint operator on $M$.

In this note we show:
\begin{theorem}
\label{main_tg}
If a closed hyperbolic $M^3$ bounds geometrically, then $\eta(M) \in {\bf Z}$.
\end{theorem}
It has been shown in \cite{MN}, that as we run over surgeries
on a hyperbolic knot in $S^3$ the $\eta$--invariant takes on a dense
set of values in $\bf R$, so we deduce:
\begin{corollary}
\label{nobound_tg}
There are closed hyperbolic 3--manifolds which do not bound geometrically.
\end{corollary}
These are the first known examples of hyperbolic $3$--manifolds,
that while they are nullbordant, are not nullbordant in this
geometrical sense; see question $4$ of \cite{W}. 

In the flat case, we show a similar theorem:
\begin{theorem}
\label{main_flat}
 If a closed flat $M^3$ is the cusp cross-section of a complete finite
 volume one--cusped hyperbolic $W^4$, then $\eta(M) \in {\bf Z}$.
\end{theorem}
In \cite{N} (see also \cite{LR}) it is shown that the homemorphism
type of every flat 3--manifold appears as a cusp cross-section of a
complete finite volume cusped hyperbolic $4$--manifold, possibly with
several cusps.  Using some calculations of \cite{Ou}, one can show
that there are flat $3$--manifolds with nonintegral $\eta$--invariant
(see the calculations below) so that we deduce:
\begin{corollary}
\label{nobound_flat}
There are closed flat 3--manifolds which are not homeomorphic to the cusp 
cross-section of any complete finite volume one--cusped  hyperbolic $4$--manifold.
\end{corollary}

The proof of Theorems \ref{main_tg} and \ref{main_flat} use a
celebrated formula of Atiyah--Patodi--Singer, but less refined
geometrical considerations still give information:
\begin{theorem}
\label{main2}
Let $M(n)$ be an infinite sequence of distinct closed hyperbolic
$3$--manifolds.  Suppose that, for each $n$, $M(n)$ is the totally
geodesic boundary of a hyperbolic $4$--manifold $W(n)$.

Then $\chi(W(n)) \rightarrow \infty$ as $n \rightarrow \infty$, where
$\chi$ denotes Euler characteristic.
\end{theorem}
Using the above and a result of Gordon \cite{Gor} we deduce:
\begin{corollary}
\label{acyclic}
There are hyperbolic integral homology 3--spheres that bound rationally
acyclic (in fact contractible) $4$--manifolds but cannot geometrically
bound any rationally acyclic hyperbolic manifold.
\end{corollary}

\section{Proofs}
The starting point of this work is the following theorem of Atiyah,
Patodi and Singer:
\begin{theorem}[See \cite{APS}, \cite{APS1}]\label{aps1}Let $W$ be a 
compact oriented Riemannian
$4$--manifold with boundary $M$ and assume that near $M$, the metric is
isometric to a product.

Then
$$ \sign(W) = \frac{1}{3}\int_W p_1 -\eta(M).$$
\end{theorem}
We briefly explain the terminology of the theorem. The left hand side
is the signature of the nondegenerate symmetric form on the
image of $H^2(W,M;{\bf Z})$ in $H^2(W;{\bf Z})$
induced via the cup product, $p_1$ is the differential $4$--form
representing the first Pontryagin class; in \cite{APS1} this is
defined as $(2\pi)^{-2}tr(R \wedge R)$, where $R$ is the curvature
matrix. 

We wish to apply this theorem in the contexts provided by \ref{main_tg}
and \ref{main_flat} and there is a technical point that in neither case
is the hyperbolic metric going to be a product near the boundary. 
However, as pointed out in \cite{APS1} page 61, (see also \cite{Gil})
the correction term necessary
if the metric is not a product near the boundary is expressible in terms of 
an integral involving the curvature and the second fundamental form $\theta$; 
this expression is given explicitly in \cite{EGH} pages 348--9 as
$$ \int_{\partial W} tr(\theta \wedge R). $$

We do not need to appeal to the details of this formula, for the two
cases that interest us we will argue that the correction term must be
zero. \par\medskip 
We now complete the proof of Theorem
\ref{main_tg}.  Suppose that $M$ is hyperbolic and bounds geometrically
the hyperbolic 4--manifold $W$.  Since we are assuming that the
boundary is totally geodesic, it follows that the second fundamental
form is zero (see \cite{BC} Theorem 2 page 194) 
and the formula of \ref{aps1} holds without
correction. Now hyperbolic manifolds are conformally flat, (one simply
notes that the standard injection $\SO_0(n,1) \rightarrow
\SO_0(n+1,1)$ preserving a codimension one totally geodesic subspace
of ${\bf H}^{n+1}$ gives a conformal action on $S^n$) and since the
Pontryagin form is a conformal invariant (see \cite{APS1}), 
it must be identically
zero. Then the formula reduces to
$$\sign(W) = -\eta(M) $$ so that $\eta(M)$ is an integer as
required.\qed \par\medskip
The proof of Theorem \ref{main_flat} is similar, but we need to argue a little
differently to compute the correction since in this case the second fundamental
form is not zero.  

Truncate the manifold $W$ with a small horoball, of Euclidean height
$k$ say.  This gives a compact 4--manifold which we denote by
$W(k)$. This manifold has flat boundary and we may write the Index formula
in this case as:
$$ \sign(W) - \frac{1}{3}\int_W p_1 +\eta(M(k)) = \correction (k),$$
where $M(k)$ is isometric to the flat manifold $M$.  As above there is
no contribution from the integral term on the left hand side.
Moreover, the $\eta$--invariant of a flat 3--manifold is independent
of choice of flat metric (\cite{Ou} and see the calculation of
$\eta(M)$ below) so that the left hand side does not depend on $k$.

Now from above, the term $\correction (k)$ is formed by integrating
over the boundary a continuous locally computable quantity which is
isometry invariant, so that the correction is bounded above by $C\cdot
{\rm volume}(M(k))$ where $C$ is a constant independent of $k$.

Since ${\rm volume}(M(k))$ tends to $0$ as the horoballs get smaller,
it follows that the correction term must be zero.

Thus $\sign(W) = -\eta(M) $ is integral as required.
\qed\par\medskip

{\bf Remarks}\qua (1)\qua In the closed case, Theorem
\ref{aps1} (which reduces to the Hirzebruch signature formula in this
case) and the same argument shows $\sign(W) = 0$ for a closed
hyperbolic $4k$--manifold.\par\medskip 

(2)\qua All known examples of
hyperbolic 4--manifolds of finite volume contain immersed totally
geodesic hyperbolic 3--manifolds which give embedded totally geodesic
hyperbolic 3--manifolds in finite covers. However, until recently, no
single example of a hyperbolic 3--manifold that {\em did} bound
geometrically was known; the first example was given
\cite{RT}. \par\medskip 

(3)\qua The smallest known hyperbolic
3--manifold with $\eta (M) \in {\bf Z}$ is the manifold identified as
$\Vol 3$, with volume that of the regular ideal simplex in ${\bf H}^3$; 
one description of this manifold being $(3,-2)$, $(6,-1)$
surgery in the Whitehead link complement (see \cite{MN}). In fact 
$\eta (\Vol 3) = 0$ (see \cite{MN}). 
Note the volume of the 3--manifold constructed 
in \cite{RT} that does bound geometrically is of the order of $200$. It
has $\eta$--invariant $0$.
\par\medskip
{\bf Proof of Theorem \ref{main2}}\qua In
even dimensions, hyperbolic volumes are basically the same as Euler
characteristics \cite{G}. In dimension $4$, the exact statement is
that if $M$ is a complete hyperbolic 4--manifold of finite volume, then
$$\Vol (M) = {4\pi^2\over 3} \chi (M).$$

Furthermore by results of Wang \cite{Wa}, there at most a finite number
of isometry classes of finite volume hyperbolic 4--manifolds of  given volume.
Thus for a given constant $C$, there are only a finite number of distinct
closed hyperbolic $4$--manifolds whose Euler characteristic is $< C$.

In the notation of \ref{main2}, suppose that there were an infinite
sequence with bounded Euler characteristics and let $D(W(n))$ be the
closed hyperbolic $4$--manifolds formed by doubling along the
boundary.  These have bounded Euler characteristics so we may pass to
a subsequence so that the doubled manifolds are all homeomorphic,
hence isometric by Mostow rigidity.

This is already a contradiction to a theorem of Basmajian \cite{B} if
the volumes of the the $M(n)$ are unbounded. If the volumes of the
$M(n)$ remain bounded, then since they are all isometrically embedded
in a single closed hyperbolic $4$--manifold, there is a global lower
bound on the injectivity radius in the sequence and this is also a
contradiction, since there are only a finite number of distinct closed
hyperbolic $3$--manifolds with injectivity radius bounded below and
volume bounded above, \cite{Th}. \qed \par\medskip
    It is observed in \cite{Gor} that $1/n$ surgery on a slice knot bounds
a contractible manifold; \ref{main2} shows that only finitely many
such manifolds could geometrically bound anything rationally 
acyclic.\par\medskip
\noindent
{\bf Calculation of $\eta (M)$ for flat 3--manifolds}\par\medskip
In \cite{Ou}, formulae for the $\eta$--invariants of Seifert fibered
3--manifolds were developed. In particular for a flat 3--manifold $M$, it is
shown that $\eta (M)$ depends only on the topology of $M$ and
is independent of the flat metric. From the theory of flat Seifert fibered
manifolds, there is a unique orientable flat 3--manifold with base for
the Seifert fibration $S^2$ and Seifert invariants $(2,1)$, $(3,-1)$, $(6,-1)$.

From \cite{Ou}, $\eta (M)$ is given by:
$$\eta (M) = 4 ({1\over 8}\cot^2 \pi/2 + {1\over {12}}
(-\cot^2 \pi/3 + \cot^2 2\pi/3)\hbox to 1.2in{\hss}$$

\vglue -20pt 

$$\hbox to 1.2in{\hss}+{1\over {24}} \Sigma_{k=1}^{5} \cot (-k\pi/6) \cot (k\pi/6).$$

This calculation gives $\eta (M) = -4/3$, as required in Corollary
\ref{nobound_flat}.

On repeating this calculation for the six other orientable flat
3--manifolds, the $\eta$--invariant is integral except for the unique
Seifert manifold with base $S^2$ and Seifert invariants $(3,2)$,
$(3,-1)$, $(3,-1)$ for which the answer is $-2/3$.  \par\medskip {\bf
Higher dimensional manifolds}\par\medskip Although the focus of this
note has been motivated by the considerations coming from
$4$--manifolds, many of our results continue to hold to obstruct
$4n-1$ manifolds as the totally geodesic or cusp cross-sections of
$4n$--manifolds. The Atiyah--Patodi--Singer theorem holds verbatim as
does the formula for the correction term. It follows that geometric
bounding in the totally geodesic case is obstructed exactly as above
by the $\eta$--invariant.  However, in contrast to the low dimensional
case, we do not know how to compute any higher dimensional
$\eta$--invariants. There is no analogue of the results of \cite{MN}
in dimensions $\geq 4$.

The case of a cusp  cross-section also works similarly,
but we need to make the additional observation (see
\cite{APS1}) that scaling the metric does not change
the $\eta$--invariant. Since the flat metrics on horoball
cross-sections at different heights differ only by
scaling, the term  $\correction(k)$ continues to be independent of $k$
and the geometrical argument above still works to show
$\correction(k) = 0$. We do not have explicit flat $4n-1$--manifolds $(n > 1)$
for which the $\eta$--invariant is non-integral, but it seems
likely that such do exist.\par\medskip
{\bf Acknowledgements}\qua The authors are indebted to Xianhze Dai for
several enormously helpful and stimulating conversations. Discussions
with C\,McA Gordon, and email correspondence with W Neumann
were also very informative.

The first author was partially supported during this work by the
NSF. The second author was partially supported by the NSF, the Alfred
P Sloan Foundation and a grant from the Texas Advanced Research
Program.

\end{document}